\documentclass[44pt]{article}
\usepackage{amssymb}
\usepackage{latexsym}
\usepackage{exscale}

\begin{document}

\def\theequation{\thesection.\arabic{equation}}
\renewcommand{\rq}[1]{(\ref{#1})}
\newcommand{\ov}{\overline}
\renewcommand{\v}{ {\,\bf v }}
\newcommand{\w}{ {\,\bf w }}
\newcommand{\dist}{\mbox{dist}}
\newcommand{\dsize}{\displaystyle}
\newcommand{\supp}{\mbox{supp}}
\renewcommand{\Re}{\mbox{Re}}
\renewcommand{\Im}{\mbox{Im}}
\renewcommand{\supp}{\mbox{supp}}

\newcommand{\C}{ \bf C }

\newcommand{\la}{\mbox{$\lambda$}}
\newcommand{\smooth}{\mbox{$\Psi^{-\infty}$}}
\newcommand{\ra}{\mbox{$\rightarrow $}}
\newcommand{\Co}{C_0^{\infty}}
\newcommand{\pa }{\partial }
\newcommand{\Da}{\Delta_{a,r}}
\newcommand{\og}{{\overline g}}
\newcommand{\Pf}{\noindent {\it Proof:}\ }
\newcommand{\ep}{\epsilon}
\newcommand{\frf}{(\frac{f'}{f})}
\newcommand{\tb}{\tilde{B}}
\newcommand{\ty}{\tilde{y}}
\newcommand{\R}{\bf R}
\newcommand{\al}{\alpha }
\newcommand{\be}{\beta }
\newcommand{\ga}{\gamma }
\newcommand{\tg}{\tilde{\gamma}}
\newcommand{\dxdt}{\frac{dX}{dT}}
\newcommand{\oa}{\overline{a_m}}
\newcommand{\bR}{{ \mathbb R  }}
\newcommand{\bN}{{ \mathbb N  }}
\newcommand{\mR}{{\mathcal R}}
\newcommand{\mS}{{\mathcal S}}
\newcommand{\mF}{{\mathcal F}}
\newcommand{\bu}{\bar{u}}
\newcommand{\f}{\footnote}
\newcommand{\BEL}{\begin{equation}\label}
\newcommand{\EE}{\end{equation}}
\newcommand {\eps}{\epsilon}
\renewcommand{\t}{\tilde}
\def \beq{\begin{equation}}
\def \eeq{\end{equation}}
\def \lab{\label}
\renewcommand{\medskip}{\vskip .5 cm}
\newtheorem{Thm}{Theorem}[section]
\newtheorem{Lemma}[Thm]{Lemma}
\newtheorem{Cor}[Thm]{Corollary}
\newtheorem{Prop}[Thm]{Proposition}
\newtheorem{Remark}[Thm]{Remark}

\title{Regularity of solutions of quasi-linear elliptic
 equations with $L\log^m L$ coefficients}

 \maketitle
\

\noindent{ to be published in Journal of Differential Equations}

\noindent{\bf Julian Edward}, Department of Mathematics and Statistics, Florida International University, Miami,
FL 33199, U.S.A.; edwardj@fiu.edu; (305)-348-3050.

 \noindent{\bf Steve Hudson},  Department of Mathematics and Statistics, Florida International University, Miami,
FL 33199, U.S.A.; hudsons@fiu.edu; (305)-348-3231.

 \noindent{\bf Mark Leckband},  Department of Mathematics and Statistics, Florida International University, Miami,
FL 33199, U.S.A.; leckband@fiu.edu .

\vspace{5mm}

\

{\bf Abstract:} Let $D$ be an bounded region
in $\bR^n$. The regularity of solutions of a family of  quasilinear elliptic partial differential equations is studied, one example being $\Delta_nu=Vu^{n-1}$. The coefficients are assumed to be in the space $L\log^{m}L(D)$ for $m>n-1$.
Using a Moser iteration argument coupled with the Moser-Trudinger inequality,  a local $L^{\infty}$ bound on the solution $u$ is proven. A Harnack-type inequality is then proven. These results are shown to be sharp with respect to $m$.
Then essential continuity of $u$ is proven, and away from the boundary
 a bound on the modulus of continuity.

\vspace{5mm}
\begin{section}{Introduction}
Let $ u\in W^{1,n}(D)$ with $D$ a bounded, connected, open subset of $ \bR^n$. Assume $u$ satisfies
the quasi-linear partial differential equation
\begin{eqnarray}
{\rm div} {\cal A} (x,u,u_x) = {\cal B} (x,u,u_x).
\label{SerrinPde}
\end{eqnarray}
on $D$. Here
$u_x=(u_{x_1},...,u_{x_n})$ is the gradient of $u$ with respect to standard coordinates $(x_1,...,x_n)$. We assume ${\cal A}(x,u,p),\ {\cal B}(x,u,p)$ are defined for all $x\in D$,  $u \in \bR$, and $p \in \bR^n$, and
\begin{eqnarray}
|\cal A| & \le & a|p|^{n-1}+ {b}|{u}|^{n-1} +e\cr
|\cal B| & \le & c|p|^{n-1}+ {d}|{u}|^{n-1} +f \cr
p\cdot \cal A & \ge & |p|^{n}- {d}|{u}|^{n}-g
\label{Serrin6}
\end{eqnarray}
where the coefficients $b$ through $g$ are nonnegative functions of $x$ and $a$ is a positive constant.
For several applications of  quasilinear pdes with this structure, the  reader is referred to (\cite{GT}, p. 260-263.) Serrin in \cite{S} studied such equations as a borderline case of a more general setting. He
assumed  $b,e\in L^{n/(n-1-\ep)}$, $c\in L^{n/(1-\ep )}$, and $d,e,f\in L^{n/(n-\ep )}$, and proved a number of results about $u$, including an $L^{\infty}$ bound, a Harnack inequality, and H\"older continuity, generalizing work on the linear case by Moser \cite{Mo}. See also \cite{Tr}.
For  more recent related works, see \cite{IKR}, \cite{NU}, \cite{SSSZ}.
The contribution of this paper is to extend Serrin's
results to   larger  function spaces. We also show that our function spaces are the largest possible
in a way that is clarified in Remark 1.3.

We define $L\log^{(n-1)/r}L(D)$ as the linear space of functions $V$ for which $ \int _D |V|\log^{(n-1)/r}(|V|+1)\;dx < \infty$. An associated norm, defined in Section 2, will be denoted $||V||_{(N_r),D}$ or just $||V||_{(N_r)}$.
A dual space with exponential norm
$||u||_{E_{nr}}$ is also defined there. Denote by $||u||_{n}$ the standard norm for $L^n$.
We will  assume
\beq
a>0 \mbox{ is constant, with }b^{n/(n-1)}, d, e^{n/(n-1)}, f, g, c^n \in L\log L^{(n-1)/r}(D)
\mbox{ and } 0<r<1.
\lab{a-g}
\eeq
Throughout this paper we will assume that
$$r>(n-1)/n.$$
 This avoids some technical issues in the proofs and simplifies the statements of our results.
 For more on this, see Sec. 2, especially the remarks following \rq{Stevenorm}.
In the case that \rq{a-g} holds for $r_1\in (0, (n-1)/n]$, it is easy to see \rq{a-g} also holds for any $r_2\in ((n-1)/n,1)$.
So, the results below still hold, but with $r_1$ replaced by $r_2 = (n-1/2)/n$, for example.
   Let
\beq \label{k}
k=(||e^{n/(n-1)}||^{(n-1)/n}_{(N_r),D}+||f||_{(N_r),D})^{1/(n-1)}+||g||_{(N_r),D}^{1/n}.
\eeq
In what follows, we denote by $B_R \subset D$ a ball of radius $R$ centered at a fixed point $z$.
When the domain of $u$ is restricted to $B_R$ we may write simply $||u||_{E_{nr},R}$ for $||u||_{E_{nr},B_R(z)}$. Unless stated otherwise, $C$ denotes a generic large constant.  It does not depend on variables such as $u$, but may depend monotonically on terms such as $||b||$.
At times, we may use $K$ or $C_1$ or $C(||b||)$, etc,  for such constants.
We may also introduce generic small constants, with similar conventions.

The first of the  following theorems appears later, slightly improved, as Theorem \ref{bounded}.
\begin{Thm}\label{bounded1}
Assume $R\le 1$.
Let $u$ be a weak solution of \rq{SerrinPde} in $D$ , with \rq{Serrin6}, \rq{a-g} and \rq{k},
and $\overline{B_{2R}}\subset D$. Then
\begin{eqnarray}
||u||_{\infty, R}\le C (||u||_{E_{nr},2R}+k),
\label{bound*}
\end{eqnarray}
\begin{eqnarray}
||u_x||_{n, R} \le C (||u||_{E_{nr},2R}+k).
\label{1uxbound}
\end{eqnarray}
\end{Thm}
For our definition of a weak solution of \rq{SerrinPde}, see Section 2.

In what follows, we use $\sup$ to denote essential supremum, and similarly for $\inf$.
\begin{Thm}\label{Harnackthm1}
Let $u\geq 0$ be a weak solution of \rq{SerrinPde} in $D$, with \rq{Serrin6} and \rq{a-g} holding. Then there exist positive constants $R_0$ and
$C$ such that if $R\leq R_0$ and $\overline{B_{8R}(z)}\subset D$, then
\begin{eqnarray}
\sup_{B_R(z)} \ u  \le C (\inf_{B_R(z)} \ u + k).
\label{Harnack1}
\end{eqnarray}
Here $C$ is independent of $ R$.
\end{Thm}

\begin{Remark}\label{LlogkL} Theorems \ref{bounded1} and \ref{Harnackthm1} are sharp in the sense that both hold for
all $0<r<1$, but
 both fail for all $r\ge 1$. We present a family of counterexamples in Section 2 for the basic equation $-\Delta_n u = V u|u|^{n-2}$. We will have $a=1$, $d=V\geq 0$,  $u\ge 0$, $b=c=e=f=g=0$, and hence $k=0$.
These examples include the borderline case $r=1$, with $d = V \in L\log^{(n-1)} L$, and $r =\infty$ with $V \in L^1$.
\end{Remark}

\begin{Thm}\label{continuousa}
Assume
the hypotheses of Theorem \ref{bounded1}.   Then one can redefine $u$ on a set of measure zero so that $u$ is continuous on $D$. If $D' \subset D$ is compact, then there exist positive constants $\gamma, K, $ and $R$ such that if $x_0\in D'$  and $0< |x-x_0| < R^2$, then
\begin{equation}
|u(x)-u(x_0)|\leq K|\log ({|x-x_0|})|^{-\gamma }.
\label{ContEqn}
\end{equation}
\end{Thm}
In Section 5 we will prove a slightly more precise version of this theorem, Theorem \ref{continuous}, where
the   dependencies of $K$, $R$ and $\gamma$ will be clarified.

Serrin in \cite{S},  with his slightly more regular
coefficients,  was able to prove H\"older continuity of $u$.
It is not clear to us whether Theorem \ref{continuousa} is the best result possible, or H\"older continuity holds. These issues are further discussed at the end of Section 5.

This paper is organized as follows. In Section 2 we discuss Orlicz spaces, and state some technical lemmas used later. We also prove Remark \ref{LlogkL}. In Section 3 we prove Theorem \ref{bounded1}, and in Section 4, Theorem \ref{Harnackthm1}, and in Section 5, Theorem \ref{continuousa}.
In an appendix we prove several technical lemmas.
\end{section}

\begin{section}{Preliminaries and simplications}
\begin{subsection}{About Orlicz and Lebesgue spaces.}
We denote the $L^p$ norm on $D$ by $||u||_p=(\int_D|u|^pdx)^{1/p}.$
Recall that a function $$\tilde{M}(t)= \int_0^t m(\tau)\;d\tau,\;\;  \;\; t \geq 0,$$ is an Orlicz function if $m(t) = \tilde{M}'(t)$ is right continuous, non-decreasing, and $m(0) = 0$.
Then $\tilde{M}(t)$ is a convex function and the resulting Luxemburg norm of $u$ is
$$||u||_{(\tilde{M})} = \inf \{ \la > 0 : \int_D \tilde{M}(\frac{|u(x)|}{\la }) \;dx \leq 1 \}.$$
For example, set $m(\tau) = e^{\tau^{s/(n-1)}}-1$ with $s>0$.
Then  $M_s(t) = \int_0^t m(\tau) \;d\tau$ and its complementary Orlicz function
$$N_s(t) = \int_0^t \log^{(n-1)/s}(\tau + 1)\;d\tau$$
form an Orlicz pair appearing often in this paper.
We cite some well known general results below.
\begin{Lemma}\label{KR} (see Krasnosel'ski\u{i}  and Ruticiki\u{i}).
For any $s>0$, and measurable functions $f,g$ on $D$,
\

A) $\int_D |f(x)g(x)|\;dx \leq 2||f||_{(M_s)}||g||_{(N_s)}.$
\vspace{2mm}

B) $ ||f||_{(N_s)} \leq \int N_s(|f(x)|) \; dx + 1 $.
\end{Lemma}
We will refer to the inequality in part A as the Orlicz-H\"older inequality.
It is understood here that both sides of either inequality can be infinite.
From
\begin{eqnarray}
t/2\log^{(n-1)/s}(t/2+1) \leq
N_s(t)
\leq t\log^{(n-1)/s} (t+1),
\label{EstimN}
\end{eqnarray}
we see that  $f\in L\log L^{(n-1)/s}$ if and only if  $||f||_{(N_s)} < \infty$.
In what follows, it will be convenient to work with the following functional:
\BEL {Stevenorm}
||v||_{E_s} = \inf \{\la >0: \int_{D} \exp(|v / \la |^{s/(n-1)})-1\ dx \le 1 \}.
\EE
For $s> n-1$,
the function $t\mapsto \exp(|t|^{s/(n-1)})-1$ is an Orlicz function, and hence $||v||_{E_s}$ will be a norm.
Later in this paper, we often set $s=rn$ and then assume $r>(n-1)/n$, to get a Luxemburg norm. This is in addition to our
necessary assumption that $r<1$.

It will often be convenient to restrict functions to a ball $B_h(x_0)$, where we may set $x_0=0$ and $B_h =B_h(0)$ for simplicity.
In this case, we denote the associated norms as
$||u||_{E_{s},h}$, etc, sometimes suppressing the $h$ if understood.
Note that $M_s(t)=\int_0^t e^{\tau^{s/(n-1)}}-1\;d\tau \leq t(e^{t^{s/(n-1)}}-1) \leq e^{C t^{s/(n-1)}}-1$.
So,
\BEL {Mineq}
||w||_{(M_s),h} \leq \inf\{k:\int_{B_h(0)} e^{C(\frac{|w| }{k})^{s/(n-1)}} - 1 \; dx \leq 1\} = C^{(n-1)/s}||w||_{E_{s},h}.
\EE
If $h\le 1$ is fixed, then by Lemma \ref{KR} and  \rq{Mineq} with $s=r$,
\BEL {Markineq}
 ||v||^n_{n,h} \leq 2||v^n||_{(M_s),h}\ ||\chi_{B_h}||_{(N_s),h} \leq C||v||^n_{E_{nr},h},
\EE
where $||\chi_{B_h}||_{(N_s),h} = (1/N_s^{-1}(1/|B_h|)$.
Since, for any small $\ep >0$, we have $N_s^{-1}(t)\geq Ct^{1-\ep /n}$, it follows that for small $h$
we have $||\chi_{B_h}||_{(N_s),h} \leq Ch^{n-\ep}$.

\vspace{5mm}

\begin{Lemma}\label{MTEn}
Let $u \in W^{1,n}_0 (B_R)$.
\vspace{2mm}

A) For $0< r <1$, $||u||_{E_{nr}} \leq C\log^{-\delta}(1+|B_R|^{-1/(1-r)})||u_x||_n,$
where $\delta = \frac{(n-1)(1-r)}{nr}$ and $C$ does not depend on $R$.

\vspace{2mm}

B) $||u||_{E_n} \leq (C_n)^{-1}(|B_R| + 1)^{(n-1)/n}||u_x||_n,$ where $C_n = (\sigma n^{(n-1)})^{1/n}$.
\end{Lemma}
These variations of the Moser-Trudinger inequality
are proved in the appendix, Section \ref{MT1}.
\vspace{5mm}

\end{subsection}
\begin{subsection}{Two technical lemmas}

For the reader's convenience, we state a lemma found in Serrin \cite{S} that we often use:
\begin{Lemma}\label{SL2}
Let $\delta>0$ and let $\alpha_i,\beta_i$ for $i=1,...,N$ be two sets of real numbers such that
$\alpha_i >0$ and $0\leq \beta_i<\delta$. Suppose that $z>0$ satisfies the inequality
$$z^{\delta}\leq \sum \alpha_iz^{\beta_i}.$$
Then
$$z\leq C  \sum \alpha_i^{\gamma_i},$$
where $C$ depends only on $N$, $\delta$, and $\beta_i$, and where $\gamma_i=(\delta -\beta_i)^{-1}$.
\end{Lemma}
\noindent
{\it \Pf:} Let max $\alpha_iz^{\beta_i} = \alpha_Iz^{\beta_I}$.
Then $z^{\delta}\leq N \alpha_Iz^{\beta_I}$, $z^{\delta-\beta_I}\leq N \alpha_I$,
$z\leq C \alpha_I^{(\delta-\beta_I)^{-1}} \le C \sum \alpha_i^{\gamma_i}.$
\hfill $\Box$

\
\vspace{2mm}

\begin{Lemma}\label{WkTestFn}
Let $ u \in W^{1,n} (D)$. Assume ${\cal A},{\cal B}$ satisfy \rq{Serrin6} and \rq{a-g}.
Suppose that for all $\phi \in C^{\infty}_0 (D)$,
\begin{eqnarray}
\int_{D} \phi_x \cdot {\cal A}(x,u,u_x) - \phi {\cal B}(x,u,u_x) \; dx = 0
\label{IntABphi}
\end{eqnarray}
Then \rq{IntABphi} also holds for all $v= \phi \in W_0^{1,n} (B_R(z))$, where $\overline{B_R(z)}\subset D$.
\end{Lemma}
The proof of a slightly more general result appears in the Appendix.
It also shows that the functions $v_x \cdot {\cal A}(x,u,u_x) $ and  $v {\cal B}(x,u,u_x) $ are integrable.
We will say $u$ is a weak solution to \rq{SerrinPde} if $u \in W^{1,n} (D)$
and \rq{IntABphi} holds.
Based on this lemma, we may refer to any $ v \in W_0^{1,n} (B_R(z))$ as an admissible test function.
\end{subsection}

\begin{subsection}{Proof of Remark \ref{LlogkL}}
To state the remark more precisely,
let $n\geq 2$. We show that for all sufficiently small $0 < \ep <1/8$,  there are positive solutions $u = u_\ep \in W^{1,n}_{0}(B_1)$ of $-\Delta_n u = V_{\ep} u^{n-1}$, such that $V_{\ep}\in L\log^{(n-1)/r}L(B_1)$ is non-negative,
with $||V_\ep||$ remaining bounded as $\ep\to 0$. Since $k=0$, the non-decreasing constants $C=C(||V_{\ep}||)$ in \rq{bound*} and \rq{Harnack1} should also remain bounded, but

\vskip 5pt
\noindent
A) (For Theorem \ref{bounded1}): ${||u||_{\infty, 1/2} \over ||u||_{{E_{rn}},1}} \to +\infty$, for $r \in [1,\infty )$, and
\vskip 5pt
\noindent
B) (For Theorem \ref{Harnackthm1}):  ${{\sup_{B_{1/8}}u} \over {\inf_{B_{1/8}}u}} \to +\infty$, for $r \in [1,\infty ]$.
\vskip 5pt
\noindent
{\it Proof:} This family of examples appears in \cite{EHL} for a different purpose.
The solutions $u=u_\ep$ are radial and positive. Let $\rho =|x|$. Then
$$ \Delta_n u(\rho ) =  (n-1)|u_{\rho}|^{n-2} ( u_{\rho \rho} + \frac{1}{\rho} u_{\rho}).$$
Let  $$ u(\rho)=\cases{a-b\rho^{\frac{n}{n-1}}
   & if $0 \leq \rho < \epsilon $ ,\cr
  -\log (\rho)  & if $\epsilon \leq  \rho \leq 1$  \cr}$$
where $a$ and $b$ are chosen below so that $u$ is differentiable.
Note that $\Delta_n u(\rho)=0$ for $\epsilon \leq  \rho \leq 1$.
Continuity at $ \rho = \epsilon$ of $u_{\rho}$ requires $ b =
\frac{n-1}{n} \epsilon^{-\frac{n}{n-1}}$, and of $u$ requires $a =
\frac{n-1}{n} - \log(\epsilon )$. We define $V_{\eps}$ by the
equation $-\Delta_n u = V_{\eps} u^{n-1} $, which gives $V_{\eps}= 0$ for  $\epsilon \leq  \rho \leq 1$, and $V_{\eps}
   \leq \frac{C\epsilon^{-n}}{u^{n-1}(\epsilon)}$  for $0 \leq \rho <
   \epsilon$. By \rq{EstimN}, $\int_{B_1(0)} N_r(V_{\eps})\;dx \leq C,$
a constant independent of $\eps$ and $r\ge 1$.
It follows by Lemma \ref{KR} part B that $||V_{\eps}||_{N_r}< C+1$ for all $\ep$. Likewise, for the case $r=\infty$, $||V_{\eps}||_{L^1}< C$.
\vskip 5pt

For part A, notice that $||u||_\infty = a \ge C|\log \epsilon|$.
Let $\la = \alpha^{-1}(|\log \epsilon |)^{\gamma}$, where $\gamma$,
$\alpha \in (0,1)$ are constants independent of $\eps$, to be specified below. We will show that
$$Q_\rho = \int_{B_\rho} \exp ((u(x)/\la )^{rn/(n-1)}) -1 \ dx \le 1$$
when $\rho=1$. Then, by \rq{Stevenorm}, $||u||_{E_{rn},1} < \la $, which implies part A.
Let $\bar{r} = rn/(n-1)$ and let $1-\gamma = 1/\bar{r}$.
For sufficiently small $\ep$ and $\al$, we have $u\le 2|\log \ep| \le |\log \ep|/\al$, and
$Q_\eps \le C\exp (|\log \epsilon |^{(1-\gamma)\bar{r}}) \epsilon^{n} = C \epsilon^{n-1} \le 1/4$.
Now, let $|x| = \rho>\epsilon$ and note that
$$(u(x)/\la )^{rn/(n-1)} = \al^{\bar{r}} |\log \rho |^{\bar{r}} \cdot |\log \ep |^{-\gamma\bar{r}}\le |\log \rho^{\alpha}|.$$
Let  $\sigma$ be the measure of the unit sphere in $\bR^{n}$. For $\al$ sufficiently small,
$$Q_1 - Q_\eps \le \sigma \int_\epsilon^1 [\exp (|\log \rho^{\alpha} |) -1] \rho^{n-1} \ d\rho \le \sigma \int_0^1 [\rho^{-\alpha}\ -1] \ d\rho =
\sigma (\frac{1}{1-\al}-1)\leq \frac{1}{2}.$$
\noindent
So, $Q_1<3/4$, proving part A.
\medskip
For part B, let $\ep <1/8$, so that $\inf_{B_{1/8}}u=\log (8)$. Since $u(0)=\frac{n-1}{n}-\log (\ep )$, part B follows.
\hfill $\Box$
\end{subsection}
\end{section}

\begin{section}{ $L^{\infty}$ bound on solution.}
In this section we prove a slightly stronger version of  Theorem \ref{bounded1}. Following Serrin \cite{S}, we first make a simplifying substitution. Let $k$ be as in \rq{k}. Set $$\bar{u}=|u|+k,\ \ \bar{b} = b+k^{1-n}e,\ \ \bar{d} = d+k^{1-n}f+ k^{-n}g.$$
We have
\begin{eqnarray}
|\cal A| & \le & a|p|^{n-1}+ \bar{b}|\bar{u}|^{n-1}, \cr
|\cal B| & \le & c|p|^{n-1}+ \bar{d}|\bar{u}|^{n-1} , \cr
p\cdot \cal A & \ge & |p|^{n}- \bar{d}|\bar{u}|^{n}.
\label{Serrin13}
\end{eqnarray}

\vspace{5mm}
\begin{Thm}\label{bounded}
Let $u$ be a weak solution of \rq{SerrinPde} in  $D$ with $\bar{b}^{n/(n-1)},c^n,\bar{d} \in L\log L^{(n-1)/r}(D)$
satisfying \rq{Serrin13}.
If $R\le 1$ and
$\overline{B_{2R}(0)}\subset D$ then,
\begin{eqnarray}
||u||_{\infty, R}\le C R^{-1/(1-r)} (||u||_{E_{nr},2R}+k),
\label{bound}
\end{eqnarray}
\begin{eqnarray}
||u_x||_{n, R} \le CR^{-1} (||u||_{E_{nr},2R}+k).
\label{uxbound}
\end{eqnarray}
\end{Thm}
\vspace{5mm}
Here $C$ is a large constant independent of $u$, $D$ and $R$.
As mentioned earlier, a similar theorem holds when $0<r\le {n-1\over n}$, but this is left to the reader.
If $u$ has compact support, Lemma \ref{MTEn} part A gives an immediate corollary, that $||u||_{\infty, R}\le K (||u_x||_n + k)$,
with $K$ independent of $u$.
\vspace{5mm}

\noindent
{\it Proof of Theorem \ref{bounded}:}
Let $R \leq h'<h \leq 2R$, and let $\eta$ be a smooth function with
$\eta (x)\in [0,1]$, such that:
\begin{equation}\label{eta}
\eta (x)=\left \{
\begin{array}{cc}
1, & |x|<h',\\
0, & |x|>h.
\end{array}
\right .
\end{equation}
Furthermore,  $\eta$ can be chosen so $|\eta_x|\leq 2(h-h')^{-1}$.
Let $l>k$ and $q\ge 1$, to be specified later; our $C$ will not depend on these.
Define a $C^1$ function $F$ by
$$F(\bar{u})=\left \{
\begin{array}{cc}
\bar{u}^q, & k\le \bar{u}\leq l,\\
ql^{q-1}\bar{u}-(q-1)l^q,& \bar{u}\geq l.
\end{array}
\right .
$$
Set $v(x)=F(\bar{u}(x)) \in W^{1,n}(D)$. Unless specified, integrals and norms will be taken over $B_h$.
\begin{Lemma}\label{FG} Let $F, v,\eta $ be defined as above. Then we have
 \begin{eqnarray}
||\eta v_x||_n^n & \leq & na \int |\eta_xv||\eta v_x|^{n -1}dx +  n q^{n-1} \int \bar{b} |\eta v|^{n-1}|\eta_x v|\ dx \cr
& + & \int c \eta v|\eta v_x|^{n-1}\ dx + C q^{n}\int \bar{d}(\eta v)^{n }\ dx
\label{serrin16}
\end{eqnarray}
\end{Lemma}
The proof integrates
a suitable test function against \rq{SerrinPde}, and uses \rq{Serrin13}.
It is given in Serrin \cite{S}, in the argument leading to his Equation 16, and the reader is referred to that paper for details.

\medskip
Returning to the proof of Theorem \ref{bounded}, we abbreviate \rq{serrin16} as
$||\eta v_x||_n^n \leq I_a + I_b + I_c +I_d$ and estimate each summand.
By H\"older's inequality,
$$I_a = na \int |\eta_xv||\eta v_x|^{n -1}dx \le na ||\eta_xv||_{n}||\eta v_x||_{n}^{n -1}.$$

We next consider $I_d$. In the rest of the paper, when $r$ is fixed, we may write $N_r$ as $N$ and
$M_r$ as $M$.
By Lemma \ref{KR} and \rq{Mineq}:
\begin{eqnarray*}
|\int \bar d(\eta v)^n\ dx | & \le & 2||\bar d||_{(N)} \ ||(\eta v)^n||_{(M)}\ ,\\
& \leq & C \ ||(\eta v)^n||_{E_r}\ ,\\
& = & C \ ||\eta v||^n_{E_{rn}}\ .
\end{eqnarray*}
So
$$I_d \le  C q^n|\int \bar d (\eta v)^n\ dx | \le C q^n||\eta v||_{E_{nr}}^n\ .$$
Reasoning as we did for $I_d$,
$$I_b \le  n q^{n-1}||\bar{b}||_{(N)} ||\eta^{n-1}\eta_xv^n||_{E_{r}} \le C   q^{n-1}||\eta_x||_\infty ||v||_{E_{rn}}^n\ .$$
Next, by the standard H\"older and  then the Orlicz-H\"older inequality with \rq{Mineq},
$$I_c \le ||cv||_n \ ||(\eta v_x)^{n-1}||_{n/(n-1)} \le ||c^n||_{(N)}^{1/n}  \ ||v||_{E_{rn}} \ ||\eta v_x||^{n-1}_{n} = C ||v||_{E_{rn}} \ ||\eta v_x||^{n-1}_{n}.$$
From \rq{serrin16}, we get
\begin{eqnarray}
||\eta v_x||_n^n \leq  C\big ( ||\eta_xv||_{n}||\eta v_x||_{n}^{n -1} +
q^{n-1}||\eta_x||_\infty ||v||_{E_{rn}}^n +
||v||_{E_{rn}} ||\eta v_x||^{n-1}_{n} +
q^n||\eta v||_{E_{rn}}^n\big ).
\label{s16'}
\end{eqnarray}
Set $K = (h-h')^{-1} + q$.
Then \rq{s16'} and \rq{Markineq} yields
\begin{eqnarray}
||\eta v_x||_n^n \leq  C(
K||v||_{E_{nr}}\ ||\eta v_x||_{n}^{n-1} +
K^n ||v||_{E_{rn}}^n +
||v||_{E_{rn}} ||\eta v_x||^{n-1}_{n} +
K^n|| v||_{E_{rn}}^n).
\label{s16'''}
\end{eqnarray}
Applying  Lemma \ref{SL2}, we get
\begin{eqnarray}
||\eta v_x||_{n} \leq CK|| v||_{E_{nr}},
\label{slemma2}
\end{eqnarray}
with $C$ independent of $h,h',q$.
Using  Lemma \ref{MTEn} part B,
and the triangle inequality,
\BEL {MT}
||\eta v||_{E_n} \le C ||(\eta v)_x||_n  \leq C(||\eta_xv||_n + || \eta v_x||_n),
\EE
Since $||\eta_x v||_{n} \leq CK|| v||_{E_{nr}}$ too, \rq{MT} implies
\BEL {rHv}
||v||_{E_n,h'} \le  C K ||v||_{E_{nr},h}.
\EE

\begin{Lemma} \label{MonoConv}
For any $h'$, and any $m>n-1$,
$$\lim_{l\to \infty}||v||_{E_m,h'}=||\bar{u}^q||_{E_m,h'}.$$
\end{Lemma}

The proof is in the Appendix. Applying the lemma to \rq{rHv} we get
$||\bar{u}^q||_{E_n,h'} \le C K||\bar{u}^q||_{E_{nr},h}$, so
\BEL {rH}
||\bar{u}||_{E_{qn},h'} \le \big ( C K\big )^{1/q} ||\bar{u}||_{E_{qnr},h}\ .
\EE
We iterate \rq{rH} with $q = q_j = r^{-j}$,
$h = h_j =R +R2^{-j}$ and $h'=h_{j+1}$, so $h-h'=R2^{-j-1}$ for $j=0,1,2\ldots$. So, $K  \le C^j/R$ and
$$||\bar{u}||_{E_{nq_j},h_{j+1}} \le  C ||\bar{u}||_{E_{nr},2R}\ \prod_{k=0}^j(C^k/R)^{r^{k}} =
C^{\sum_0^{j}kr^k} R^{-\sum_0^{j}r^k} ||\bar{u}||_{E_{nr},2R} \le C R^{-1 /(1-r)}||\bar{u}||_{E_{nr},2R}$$
since $r<1$ and $R\le 1$.

\begin{Lemma} \label{boundedE} Assume
$||\bar{u}||_{E_{q_j},h_{(j+1)}} \le J$ for all $j$ and that $q_j\to \infty$.  Then $\bar{u}\in L^\infty(B_R)$ with
$$||\bar{u}||_{\infty, R} \le  J.$$
\end{Lemma}
The proof is in the Appendix.
Since $ ||\bar{u}||_{E_{nr},2R} \le ||u||_{E_{nr},2R}+ k||1||_{E_{nr},2R}$ and $||1||_{E_{nr},2R} \le C$
this proves \rq{bound}. Also, \rq{uxbound} follows from \rq{slemma2} by setting $q=1$, $h=2R$, $h'=R$
and $\eta\equiv 1$ on $B_R(0)$,
so that
$$||u_x||_{n,R} \leq CK||\bar u||_{E_{nr},2R} \le CR^{-1}(||u||_{E_{nr},2R}+k).$$
This proves Theorem \ref{bounded}.
\hfill $\Box$
\end{section}

\begin{section}{Harnack inequality}

Assume $u\geq 0$ is a weak solution of \rq{SerrinPde} in $D$, with \rq{Serrin6} and \rq{a-g} holding.
We assume $z\in D$ and $\overline{B_{8R}(z)} \subset D$,
 and
$${k}=(||{e}^{n/(n-1)}||^{(n-1)/n}_{(N_r),B_{8R}}+||{f}||_{(N_r),B_{8R}})^{1/(n-1)}+||{g}||_{(N_r),B_{8R}}^{1/n}.$$
 In this section we prove Theorem \ref{Harnackthm1}, i.e. for small enough $R$,
\begin{eqnarray}
\sup_{B_{R}(z)} \ u  \le C (\inf_{B_{R}(z)} \ u + k).
\label{Harnack}
\end{eqnarray}
Here
$C$ may depend on the parameters below, but not on $u$, $z$ or $D$:
\begin{equation}
n,r, a, ||b^{n/(n-1)}||, ||d||, ||e^{n/(n-1)}||, ||f||, ||g||, ||c^n||, \label{Cdep}
\end{equation}
all these norms being $||*||_{(N_r),D}$. The key proposition, proven in the next subsection is:
\begin{Prop}
\label{Harnackprop} Assume $a,\bar{b},c,\bar{d}$  satisfy
the hypotheses of Theorem \ref{bounded}.
Then there exists $R_0>0$ depending on the quantities given in \rq{Cdep}, such that if
$\overline{B_{8R_0}(z)} \subset D$, then \rq{Harnack} holds with $R=R_0$.
\end{Prop}
To prove Theorem \ref{Harnackthm1}, we need to show that \rq{Harnack} holds for all $R<R_0$, with $C$ independent of $R$.
The dilation argument for this  appears in Subsection [4.2]. Then, an easy finite cover argument gives:
\begin{Cor}
Let $D'$ be any compact subset of $D$. Then there exists a constant $C$ depending on $D'$, $D$ and quantities in \rq{Cdep}
such that $$\sup_{D'} \ u  \le C (\inf_{D'} \ u + k).$$
\end{Cor}

\begin{subsection}{Proof of Proposition \ref{Harnackprop}}
\
\noindent
{\it Proof of Proposition \ref{Harnackprop}.}
With $u$ as above, let $\bu = u+k+\ep$ with $\ep> 0$, so $\bu ^{-1} \in W^{1,n}(D)$.
 We can assume that $z=0$, and will denote $B_{R_0}(0)$ by $B_{R_0}$. We will specify $R_0$ in the proof of Lemma \ref{LebToExp}. Below, $C$ is a generic large constant determined by the quantities in \rq{Cdep} and by $R_0$.
The proposition follows from the following string of inequalities. We prove them below, and then let $\ep \to 0$.
\begin{eqnarray}
 ||u ||_{\infty, R_0} & \le & C||\bu||_{E_{n},2R_0} \label{x1}\\
& \le & C ||\bu ||_{n/r,4R_0} \label{x2}\\
& \le & C ||\bu^{-1} ||_{n/r,4R_0}^{-1}\label{x3}\\
& \le & C ||\bu^{-1} ||_{E_{n},2R_0}^{-1}\label{x4}\\
& \le & C ||\bu^{-1}||^{-1}_{\infty, R_0} = C \inf_{B_{R_0}}\ \bu \le C (\inf\ u + k+\epsilon )\label{x5}
\end{eqnarray}

Inequality \rq{x1} follows easily from the proof of Thm \ref{bounded}, by omitting the last few steps of its proof; for these steps  the only requirement on $R_0$ is that $\overline{B_{2R_0}} \subset D$ . To prove \rq{x2} and \rq{x4} we set $v(x)=\bu (x)^q$ with $q\not =0$, possibly negative.
Let $\beta =nq-n+1$ and let $$\phi (x)=\eta^n(\bu (x))^{\beta},$$
with $\eta$ as in \rq{eta}, except that now $h'=2R_0$ and $h=4R_0$.
Both $\bu ^{-1}$ and $\bu \in L^{\infty}\cap W^{1,n}$
(by Theorem \ref{bounded}). So, $\phi \in W_0^{1,n} (B_{h})$ is an admissible
test function for any $q$ by Lemma \ref{WkTestFn}.
The following is an analogue of \rq{serrin16}.
Unless specified, integrals and norms in this section
are over $B_h$.
\begin{Lemma}\label{serrin16'}
For any $q\neq 0$,
\begin{eqnarray*}
\frac{|\beta|}{|q|}||\eta v_x||_n^n  & \leq &  na \int |\eta_xv||\eta v_x|^{n -1}dx +  n |q|^{n-1} \int \bar{b} |\eta v|^{n-1}|\eta_x v|\ dx \\
& + & \int c |\eta v||\eta v_x|^{n-1}\ dx + (1+|\beta | ) |q|^{n -1}\int \bar{d}(\eta v)^{n }\ dx  .\\
\end{eqnarray*}
\end{Lemma}
{\it Proof of Lemma \ref{serrin16'}:}  We have
$$\phi_x\cdot {\cal A}+\phi {\cal B}=  (\eta^n\beta \bu^{\beta -1}u_x+n\eta^{n-1}\bu^{\beta}\eta_x)\cdot {\cal A}+\eta^n\bu^{\beta}{\cal B}.$$
  The following formulas hold:
\begin{eqnarray*}
\bu^{\beta}\bu^{n-1} & = &v^n,\\
|q|^{n-1}|\bu|^{\beta}|u_x|^{n-1} & = & v|v_x|^{n-1} ,\\
|q|^n|\bu|^{\beta -1}|u_x|^n& = & |v_x|^n.
\end{eqnarray*}
{\bf Case 1: $\beta >0$.}
In this case, by  \rq{Serrin13} we have
\begin{eqnarray}
\phi_x\cdot {\cal A}+\phi {\cal B}& \geq & \eta^n\beta |\bu|^{\beta -1}({|u_x|^n-\bar{d}\bu^n})-n\eta^{n-1}|\bu|^{\beta}|\eta_x|(a|u_x|^{n-1}+\bar{b}|\bu |^{n-1})-\eta^n|\bu|^{\beta}(c|u_x|^{n-1}+\bar{d}|\bu |^{n-1})\nonumber \\
& = &   \frac{|\beta|}{|q|^n}|\eta v_x|^n-\frac{an}{|q|^{n-1}}|\eta v_x|^{n-1}|\eta_x v|-n\bar{b}|\eta v|^{n-1}|\eta_xv|- \frac{c}{|q|^{n-1}}|\eta v||\eta v_x|^{n-1}-\bar{d}|\eta v|^n(1+|\beta | ).\nonumber \\
& & \label{cal1}
\end{eqnarray}
Integrating, we get the lemma in this case.

\vskip 2mm
\noindent
{\bf Case 2: $\beta <0$.}
In this case, the argument in \rq{cal1} goes as follows:
\begin{eqnarray*}
\phi_x\cdot {\cal A}+\phi {\cal B}& \leq & -\eta^n|\beta| (\bu)^{\beta-1}({|u_x|^n-
\bar{d}\bu^n})+n\eta^{n-1}\bu^{\beta}|\eta_x|(a|u_x|^{n-1}+\bar{b}|\bu |^{n-1})+\eta^n\bu^{\beta}(c|u_x|^{n-1}+\bar{d}|\bu |^{n-1})\nonumber \\
& = &   -\frac{|\beta |}{|q|^n}|\eta v_x|^n+\frac{an}{|q|^{n-1}}|\eta v_x|^{n-1}|\eta_x v|+n\bar{b}|\eta v|^{n-1}|\eta_xv|+ \frac{c}{|q|^{n-1}}|\eta v||\eta v_x|^{n-1}+\bar{d}|\eta v|^n(1+|\beta | ).
\end{eqnarray*}
Integrating again, the lemma follows.
\hfill $\Box$
\medskip
The following lemma proves inequalities \rq{x2} and \rq{x4}.
\begin{Lemma}  \label{LebToExp}
There exists $R_0>0$  such that
\begin{eqnarray}
||\bu ||_{E_{n},2R_0} \leq C ||\bu ||_{n/r,4R_0} \ \ {\it and}
\label{L2E+}
\end{eqnarray}

\begin{eqnarray}
C ||\bu^{-1} ||_{E_{n},2R_0}^{-1} \geq  ||\bu^{-1} ||_{n/r,4R_0}^{-1}\ .
\label{E2L-}
\end{eqnarray}
\end{Lemma}

\noindent
{\it Proof:} Below,  $C_1$ will denote a large generic constant depending only on $r,n,a,  ||\bar{b}^{n/(n-1)}||_{(N_r),D}, \ ||{\bar{d}}||_{(N_r),D}$ and
$||{{\bar{c}}}^n||_{(N_r),D}$. We estimate the terms on the right hand side of Lemma \ref{serrin16'}. First,
\begin{eqnarray*}
\int \bar{b} |\eta v|^{n-1}|\eta_x v|\ dx & \leq & \big ( \int \bar{b}^{n/(n-1)}|\eta v|^n \ dx\big )^{(n-1)/n} \ ||\eta_x v||_n,\\
& \leq & || \bar{b}^{n/(n-1)}||_{(N)}^{(n-1)/n}\ ||(\eta v)^n||_{E_{r}}^{(n-1)/n}\ ||\eta_x v||_n,\\
& \leq & C_1 ||\eta v||_{E_{rn}}^{n-1}\ ||\eta_x v||_n,\\
& \leq & C_1 ( ||\eta v||_{E_{rn}}^{n} + ||\eta_x v||_n^n).
\end{eqnarray*}
We used Young's Inequality in the last step. Much as in the proof of \rq{s16'},
$\int |\eta_xv||\eta v_x|^{n -1}dx \leq ||\eta v_x||_n^{n-1} ||\eta_xv||_n$, $\ \int c |\eta v||\eta v_x|^{n-1}\ dx\leq C_1 ||\eta v_x||_n^{n-1}||\eta v||_{E_{rn}}$ and
$\int \bar{d}(\eta v)^{n }\ dx\leq C_1 ||\eta v||_{E_{rn}}^n.$
By Lemma \ref{serrin16'}, regarding $|q|=1/r$ and $|\beta|$ as fixed constants,
$$
||\eta v_x||_n^n   \leq  C_1 \left (||\eta v_x||_n^{n-1}( ||\eta_xv||_n +   ||\eta v||_{E_{rn}})
 + ||{\eta} v||_{E_{rn}}^n+||\eta_x v||_n^n \right ).
$$
By Lemma \ref{SL2}, we get $||\eta v_x||_n \leq C_1 (||\eta v||_{E_{rn}}+||\eta_x v||_{n})$.
Now by  Lemma \ref{MTEn} part A,
\begin{eqnarray*}
||\eta v||_{E_{rn}} & \leq & C(h)||(\eta v)_x||_n\\
& \leq & C(h)(||\eta_x v||_n+||\eta v_x||_n) \\
& \leq & C(h)\big ( ||\eta_x v||_n+ C_1(||\eta v||_{E_{rn}}+||\eta_x v||_{n})\big ) .
\end{eqnarray*}
Again by Lemma \ref{MTEn} part A, we can choose
$R_0=h/4$ sufficiently small that $C(h) C_1<1/2$, so
$$||\eta v||_{E_{rn}} \leq 2C(h) (1+C_1) ||\eta_x v||_n\le C||\eta_x v||_n .$$
To prove \rq{L2E+}, let $q=1/r$ so that
$||\bu ||_{E_{n},2R_0}^q =
||\bu^q ||_{E_{rn},2R_0}  \le ||\eta v||_{E_{rn}}$,
and $||\eta_x v||_{n} \leq C ||\bu^q ||_{n} = C ||\bu ||_{n/r}^q$.
Likewise, \rq{E2L-} follows from $q=-1/r$.
\hfill $\Box$
\vspace{5mm}

\noindent
In the next passage, we return to the proposition, proving \rq{x3}.
The Proposition \ref{MTavg} below is based on a well-known version (see \cite{Leck}) of the Moser-Trudinger inequality: if $||w_x||_{n,B}\le 1$ then
\begin{equation}
   \int_B \exp(\beta_n|w-w_B|^{n\over n-1})\ dx \le C_n |B|,\label{MT9}
\end{equation}
where $\beta_n, C_n$ are positive constants independent of $w$, $B$ is any ball, and $w_B$ is the mean value.

\begin{Prop}\label{MTavg}

 Suppose $w\in W^{1,n}(B)$ with $||w_x||_{n} \leq T$. For every $p>0$,
 $$\int_B e^{pw/T}\ dx \int_B e^{-pw/T}\ dx \le C|B|^2,$$
where $C=C_n+e^{(Tp)^n\beta_n^{1-n}}$.
\end{Prop}

\vspace{5mm}
\noindent
{\it Proof:} Fix $p>0$. Then

 \begin{eqnarray*}\int_B e^{pw}\ dx \int_B e^{-pw}\ dx
 & = &\int_B e^{p(w-w_B)}\ dx \int_B e^{-p(w-w_B)}\ dx \\
&\le & \int_B e^{p|w-w_B|}\ dx \int_B e^{p|w-w_B|} \ dx\\
&  = &[\int_B e^{p|w-w_B|} \ dx]^2 .
\end{eqnarray*}

\vspace{5mm}
\noindent
 Partition $B = B_1\cup B_2$ with $B_1 = \{x\in B: \beta_n |w(x)-w_B|^{1\over n-1} \ge pT^{n/(n-1)}\}$ and
 $B_2 = \{x\in B: \beta_n|w(x)-w_B|^{1\over n-1} < pT^{n/(n-1)}\}$.

\vspace{5mm}
\noindent Then

$$\int_{B_1} e^{p|w-w_B|} \ dx \le \int_{B_1} e^{\beta_n |(w-w_B)/T|^{1+{1\over n-1}}} \ dx \le C_n |B|$$
by \rq{MT9}. And

$$\int_{B_2} e^{p|w-w_B|} \ dx \le \int_{B_2} e^{p (p/\beta_n )^{n-1}T^n} \ dx \le e^{T^np^n \beta_n^{1-n}}|B|  $$
This proves the proposition.
\hfill $\Box$

\


\begin{Lemma}\label{+to-}
For any $p_0>1$, there exists a constant $C$
such that
\begin{equation}
||\bu ||_{p_0,4R_0}\leq C ||(\bu )^{-1}||^{-1}_{p_0,4R_0}.
   \label{Lemma+to-}
\end{equation}

\end{Lemma}
\noindent
{\it Proof:} Define $\eta$ as in \rq{eta} but with $h=8R_0$ and $h'=4R_0$.
Set $\phi (x)=\eta^n \bu^{1-n}$, which is admissible by Lemma \ref{WkTestFn}.
The argument in \cite{S} (p.266, Case IV) leading to his Eq.36 shows that
\medskip

\BEL{Serrin36}
(n-1)||\eta v_x||_n^n \leq na\int |\eta_x||\eta v_x|^{n-1}dx\ + \ n\int\bar{b}\eta^{n-1}|\eta_x|dx\ +\  \int c\eta|\eta v_x|^{n-1}dx\ +\
n\int \bar{d}\eta^n dx,
\EE
where $v(x)=\log (\bu (x))$. Calculations similar to those in Theorem \ref{bounded} and Lemma \ref{LebToExp} show
$||\eta v_x||_{n,8R_0}^n\leq C(1+ ||\eta v_x||_{n,8R_0}^{n-1})$.
Hence by Lemma \ref{SL2}, $|| v_x||_{n,4R_0}\leq C$.
Applying Proposition \ref{MTavg}, with $w=v$, $p= p_0C $ and $T=C$,
 $$\int_{B_{4R_0}} e^{p_0v}\ dx \int_{B_{4R_0}} e^{-p_0v}\ dx \le C.$$
Since $v=\log(\bu)$,
inequality \rq{Lemma+to-} follows.
\hfill $\Box$
\vskip 10pt
Lemma \ref{+to-} with $p_0=n/r$ proves \rq{x3}.
Next, we prove \rq{x5}.
\begin{Lemma}\label{inflb} We have
\begin{eqnarray}
||(\bar{u})^{-1}||_{E_{n},2R_0}^{-1} \leq C(\inf_{B_{R_0}} u+k+\ep)  \ .
\label{Min2Ej}
\end{eqnarray}
\end{Lemma}

\noindent
{\it Proof:} Let $R_0 \leq h'<h \leq 2R_0$ and choose $\eta$ as in \rq{eta}.
In this proof, $C$ is a constant depending as usual
on the terms in \rq{Cdep}, but independent of $\bar{u}$, $h'$, $h$ and $q$.
Assume $q\leq -1/r$ and notice that $q^{-1}\beta >1$. By Lemma \ref{serrin16'},
\begin{eqnarray}
||\eta v_x||_n^n & \leq & na \int |\eta_xv||\eta v_x|^{n -1}dx +  n |q|^{n-1} \int \bar{b} |\eta v|^{n-1}|\eta_x v|\ dx \cr
& + & \int c \eta v|\eta v_x|^{n-1}\ dx + C|q|^{n}\int \bar{d}(\eta v)^{n }\ dx, \cr
& = & I_a+I_b+I_c+I_d,
\label{serrin16har}
\end{eqnarray}
with $v=(\bu)^q$.
Calculations as in the proof of Theorem 1, based on H\"older's inequality, the Orlicz-H\"older inequality,
Lemmas \ref{KR} and \rq{Mineq} show that
$I_a \le na ||\eta_xv||_{n}||\eta v_x||_{n}^{n -1},$
$I_b \le C |q|^{n-1}||\eta_x||_\infty ||v||_{E_{rn}}^n,$
$I_c \le C ||v||_{E_{rn}} \ ||\eta v_x||^{n-1}_{n}$ and
$I_d \le C |q|^n||\eta v||_{E_{nr}}^n.$
Set $K = (h-h')^{-1} + |q|$.
From \rq{serrin16har} and \rq{Markineq}

\begin{eqnarray}
||\eta v_x||_n^n \leq  C(
K||v||_{E_{nr}}\ ||\eta v_x||_{n}^{n-1} +
K^n ||v||_{E_{rn}}^n +
||v||_{E_{rn}} ||\eta v_x||^{n-1}_{n} +
K^n|| v||_{E_{rn}}^n).
\label{Hars16'''}
\end{eqnarray}
Applying  Lemma \ref{SL2}, $||\eta v_x||_{n} \leq CK|| v||_{E_{nr},h}$.
By Lemma \ref{MTEn} part B, the triangle inequality
and  \rq{Markineq}, we get
\begin{eqnarray}
||v||_{E_n,h'} \le ||\eta v||_{E_n,h}
& \leq & C ||(\eta v)_x||_{n,h}  \cr
& \leq & C(||\eta_xv||_{n,h} + || \eta v_x||_{n,h}) \cr
& \leq & C K ||v||_{E_{nr},h}.
\label{rHvb}
\end{eqnarray}
Thus
$||\bar{u}^q||_{E_n,h'} \le C K||\bar{u}^q||_{E_{nr},h}$, and we have
$$||(\bar{u})^{-1}||_{E_{n|q|},h'} \le (C K)^{1/|q|}||(\bar{u})^{-1}||_{E_{nr|q|},h}.$$
Thus
\BEL {rHb}
||(\bar{u})^{-1}||_{E_{n|q|},h'}^{-1} \ge (C K)^{-1/|q|}||(\bar{u})^{-1}||_{E_{nr|q|},h}^{-1}.
\EE
We iterate \rq{rHb} with $h = h_j =R_0 +R_02^{-j}$ and $h'=h_{j+1}$, so $h-h'=R_02^{-j-1}$ for $j=0,1,2\cdots$.
Also, let $q = q_j = - r^{-j-1}$,
and note that
$$K_j = 2^{j+1}/R_0 + (1/r)^{j+1} \le C^{j+1}.$$

Thus
\begin{eqnarray*}
||(\bar{u})^{-1}||_{E_{nr^{-k-1}},h_{k+1}}^{-1} & \ge &   \prod_{j=0}^k({ C}^{j+1})^{-r^{j+1}}\ ||(\bar{u})^{-1}||_{E_{n},2R_0}^{-1}\\
&  = & { C}^{-r(\sum_1^{k}(j+1)r^j)}\ ||(\bar{u})^{-1}||_{E_{n},2R_0}^{-1}\\
&  \ge & C_2 ||(\bar{u})^{-1}||_{E_{n},2R_0}^{-1}
\end{eqnarray*}
with a generic small constant  $C_2>0$ because $r<1$. Evidently, $C_2$ depends only on the quantities listed in \rq{Cdep} together with $R_0$. We conclude by Lemma \ref{boundedE} that
$$\inf_{B_{R_0}}\bu =||(\bar{u})^{-1}||_{{\infty},R_0}^{-1}\ge C_2 ||(\bar{u})^{-1}||_{E_{n},2R_0}^{-1}\ .$$
Since $ \bar{u} =u+ k+\epsilon$, this proves \rq{x5}
and completes the proof of Proposition \ref{Harnackprop}.
\hfill $\Box$
\medskip

\begin{Remark}
Careful reading of the proof of Proposition \ref{Harnackprop} shows that the constant $C$ appearing in \rq{Harnack} can be chosen to be a non-decreasing function in each of the arguments $||\bar{b}||,\ ||c||,\ ||\bar{d}||$.
\end{Remark}
\end{subsection}

\begin{subsection}{Proof of  Theorem \ref{Harnackthm1}}
Let $z\in D$ and $R_0$ be as in Proposition \ref{Harnackprop}, but we may assume $z=0$. Choose $R\leq R_0$ so that $\overline{B_{8R}(z)}\subset D$.
As a first step to proving Theorem \ref{Harnackthm1},
we show that our problem can be rescaled to a disk of radius $8R_0$.
Let $u$ be a weak solution of \rq{SerrinPde} in the ball $B_{8R}(0)$.
Let  $\rho =R/R_0 \le 1$ and $v(x)=u(\rho x)$.
Thus $v_x = \rho u_x(\rho x)$, and $v$ is a weak solution in the ball $B_{8R_0}(0)$ of
$$
{\rm div}\; {\cal C}(x,v_x,v) + {\cal D}(x,v_x,v) = 0,
$$
where ${\cal C} (x, p, u) = \rho^{n-1}{\cal A} (\rho x, p/\rho , u)$ and ${\cal D} (x, p, u) = \rho^{n} {\cal B} (\rho x, p/\rho , u)$. The inequalities of \rq{Serrin13} become
\begin{eqnarray}
|\cal C| & \le & a|p|^{n-1}+ \rho^{n-1}{\bar{b}(\rho x)}|{v}|^{n-1} = a|p|^{n-1}+ \hat{b}(x)|{v}|^{n-1} \cr
|\cal D| & \le & \rho {c}(\rho x)|p|^{n-1}+ \rho^n{\bar{d}(\rho x)}|{v}|^{n-1}   = \hat{c}(x)|p|^{n-1}+ \hat{d}(x)|{v}|^{n-1}  \cr
p\cdot \cal C & \ge & |p|^{n}- \rho^n{\bar{d}(\rho x)}|{v}|^{n}=|p|^{n}- \hat{d}(x)|{v}|^{n},
\label{RSerrin6*}
\end{eqnarray}
where  $\hat{b},\hat{c},$ and $\hat{d}$ have been implicitly defined.
\begin{Lemma}\label{dilate*}
Let $0 <R \leq R_0$.
Let $h$ be any one of the coefficients $\bar{b}^{n/(n-1)}, \bar{d},  {c}^n$ and let $\hat{h}$ be its counterpart in \rq{RSerrin6*}.
Then, $$||\hat{h}||_{(N_r),B_{8R_0}} \leq  ||h||_{(N_r),B_{8R}}.$$
\end{Lemma}

\noindent
{\it Proof:}
We show this for $h(x) = \bar{d}(x)$ with $\hat{h}(x) = \hat{d} = \rho^n\bar{d}(\rho x)$, noting the arguments for the other coefficients are similar. We begin with $||\hat{d}||_{(N_r),B_{8R_0}} = \inf\{\la : \int_{B_{8R_0}} N_r(\frac{\rho^n\bar{d}(\rho x)}{\la})\;dx \leq 1 \}$.
Since $N_r$ is convex, $N_r(\frac{\rho^n{\bar{d}(\rho x)}}{\la }) \leq \rho^nN_r(\frac{\bar{d}(\rho x)}{\la })$. So, we have
$$||\hat{d}||_{(N_r),B_{8R_0}} \leq  \inf\{\la : \int_{B_{8R_0}} \rho^nN_r(\frac{\bar{d}(\rho x)}{\la })\;dx \leq 1 \} = ||\bar{d}||_{(N_r),B_{8R}}.\ \ $$ \hfill $\Box$
\

We now complete the proof of the theorem. Let $u$, $v$, $R$, $R_0$ be as above. Let
$$\hat{k}=(||\hat{e}^{n/(n-1)}||^{(n-1)/n}_{(N_r),B_{8R_0}}+||\hat{f}||_{(N_r),B_{8R_0}})^{1/(n-1)}+||\hat{g}||_{(N_r),B_{8R_0}}^{1/n}.$$
By Proposition \ref{Harnackprop},
$$
\sup_{B_{R}} \ u =
\sup_{B_{R_0}} \ v  \le
C_1 \inf_{B_{R_0}} \ v+\hat{k} \le
C_1 \inf_{B_{R}} \ u+\hat{k}
$$
Here
\begin{equation}
C_1=C_1(r,n,a, ||\hat{b}^{n/(n-1)}||_{(N_r),B_{8R_0}}, \ ||\hat{d}||_{(N_r),B_{8R_0}}, \ ||\hat{c}^n||_{(N_r),B_{8R_0}}),\label{tCSteve}
\end{equation}
is the constant arising in Proposition \ref{Harnackprop} when applied to $v$.
By Remark 4.9, we can assume $C_1$ to be a non-decreasing function in its last three arguments. Hence,
by Lemma \ref{dilate*}, $C_1 \le C$, a constant that depends only on the quantities in \rq{Cdep}. Also by Lemma \ref{dilate*}, we have
  $\hat{k}\leq k$. This proves Theorem \ref{Harnackthm1}.
  \hfill $\Box$

\end{subsection}
\end{section}

\begin{section}{Continuity}
In this section, we mainly prove Theorem \ref{continuous}, a version of Theorem \ref{continuousa}. At the end of the section, we compare this theorem with its counterpart in \cite{S}, also \cite{Tr}. The next lemma shows our Orlicz spaces are nested, with an inequality similar to a standard one for Lebesgue spaces based on H\"older's inequality.
\medskip
\noindent
\begin{Lemma}\label{shrinkS}
Let
$\ep =(1-r)/2$, and  $G = G(R) = |\log(1/R)|^{-\ep}$. Suppose $||h||_{(N_r),D}<\infty$.
Then there exist $R_1<3^{-10}$ and
$C>0$, both independent of $h$ and $C$ independent of $R$, such that
$$||h||_{(N_{r+\ep}),B_R} < CG(R)||h||_{(N_r),B_R},\ \forall R<R_1.$$
\end{Lemma}
We remark that the restriction $R_1<3^{-10}$ will arise in the proof of Lemma 5.3.
\medskip
\noindent
{\it Proof:}
We will specify $R_1$ below. Fixing $R<R_1$, we have $G\leq 1$, since $R_1<1/e$. All integrals and norms here will be over $B_R$. We may assume that $h\ge 0$ and that $||h||_{(N_r)}=1$.
It will be useful to define an alternative norm
\begin{equation}
||h||_{N_s,B_R} = \inf_{\la >0}{\la  (1 + \int_{B_R} N_s(|h|/\la  )\;dx)} ,\ (n-1)/n< s<1.
\label{normNr}
\end{equation}
We have the following equivalence of norms: $||h||_{(N_s),B_R} \leq ||h||_{N_s, B_R} \leq 2||h||_{(N_s),B_R}$ (see \cite{KR}).
Applying the standard H\"older inequality,
$$\int h\log^{(n-1)/(r+\ep)}(h/G\ +1)\ dx \le [\int h\log^{(n-1)/r}(h/G\ +1)\ dx]^{r/(r+\ep)}\cdot [\int h \ dx]^{\ep/(r+\ep)} :=AB.$$
We first estimate $B$.  The H\"older-Orlicz inequality gives
\begin{equation}
\int h \;dx \leq 2||h||_{(N_{r})}||\chi_{B_R}||_{(M_{r})} = 2C\log^{-(n-1)/r}(1/|B_R|+1).\label{HO}
\end{equation}
We may assume $|B_{R_1}|\le R_1$, so
\begin{equation}
\int h \ dx \le C \log^{-(n-1)/r}(1/R),\label{hint}
\end{equation}
and $B\le C G^{(n-1)/(r(r+\ep))} \le C G$, with  $C$ depending only on $n,r$.

Next, we estimate $A$.
Factoring $h/G +1 = {h+G\over 2} \cdot {2\over G}$ and using $(X+Y)^p\leq 2^p(X^p+Y^p)$,
$$A\le C[\int h\log^{(n-1)/r}(h/2 +1)\ dx + \log^{(n-1)/r}(2/ G) \int h\ dx  ]^{r/(r+\ep)}.$$
Since $||h||_{(N_r)}=1$, we have
$\int N_r(h)\ dx =1$. By \rq{EstimN}, $h\log^{(n-1)/r}(h/2\ +1) \le 2N_r(h)$.
Hence $\int h\log^{(n-1)/r}(h/2+1)\ dx \le 2\int N_r(h)\ dx \le 2$.
For sufficiently small $R_1$, $R$ and $G$ must be small enough that, by \rq{hint}
$$\log^{(n-1)/r}(2/G) \int h\ dx \le \sqrt{G} \le 1.$$
We have
\begin{equation}\label{goal}
\int h\log^{(n-1)/(r+\ep)}(h/G\ +1)\ dx \le AB \le CG
\end{equation}
with $C$ depending only on $n,r$. We conclude by \rq{EstimN} that
$$\int N_{r+\ep}(h/G)\ dx \le C.$$
By \rq{normNr} with $\la =G$, and the remark above on norm equivalence,
$$||h||_{(N_{r+\ep})} \le ||h||_{N_{r+\ep}} \le G (1+\int N_{r+\ep}(h/G)\ dx) \le CG = CG||h||_{(N_{r})}.$$
\hfill $\Box$
\medskip

\begin{Thm}\label{continuous}
Assume the hypotheses of Theorem \ref{bounded1}.   Then one can redefine $u$ on a set of measure zero so that $u$ is continuous on $D$. If $D' \subset D$ is compact, $x_0\in D'$,  and $0< |x-x_0| < R^2$, then
\begin{equation}
|u(x)-u(x_0)|\leq K|\log ({|x-x_0|})|^{-\gamma }.
\label{ContEqn}
\end{equation}
\end{Thm}

Here, $R$ is any small positive constant such that $\overline{B_{8R}(x_0)}\subset D$ for all $x_0\in D'$, and such that
$R\le$ min $\{R_0, R_1\}$.
For $R_0$ see Theorem \ref{Harnackthm1} and for $R_1$ see Lemma \ref{shrinkS}.
Let
$$D''=\overline{\cup_{x_0\in D'}B_R(x_0)},$$
a compact subset of $D$. By Theorem \ref{bounded}, $u\in L^{\infty}(D'')$.
Then $K$ depends on $||u||_{\infty,D''},a,..., g, D, D', r, R_0$ and $R_1$.
In our proof, $\gamma = (1-r)/2n$. One could use any $\gamma < (1-r)/n$, but then $K$ may depend on $\gamma$.

\medskip
\noindent
{\it Proof:} We aim for a bound on oscillation similar to \rq{ContEqn} which does not assume continuity of $u$ (see Lemma \ref{omegabound}).
We may ignore a set of measure zero until the last paragraph of the proof, and will denote the essential supremum of $u$ by $\sup u$.
Below, all balls will be centered at $x_0$.
For all $\rho \leq R$, set
$$M(\rho )= \sup_{B_{\rho}}u,\ \mu (\rho )=\inf_{B_{\rho}}u. $$
Fixing $\rho$ for now,
$$\bu(x) =M(\rho)-u(x)\ge 0 \mbox{ and } \bar{\bu}(x)=u(x)-\mu(\rho)\ge 0$$
on $B_\rho$. Also, $\bu$ satisfies
\begin{eqnarray}
{\rm div} \bar{\cal A} (x,\bu,\bu_x) = \bar{\cal B} (x,\bu,\bu_x).
\label{SerrinPdebar}
\end{eqnarray}
with
$\bar{\cal A} (x,\bu,\bar{p})={\cal A} (x,M-\bu,-\bar{p})$ and a similar formula for $\bar{\cal B}$. These quantities obey the
inequalities
$$|\bar{\cal A} |\leq a|\bar{p}|^{n-1}+\bar{b}|\bu |^{n-1}+\bar{e}, \mbox{ etc},$$
where $\bar{b} =2^nb, \ \bar{e}=e+2^nb||u||_{\infty, D''}^{n-1}$, etc.
Let $\ep $ be as in Lemma ~\ref{shrinkS}, and
set
\begin{equation}\label{bark1}
\bar{k}=\bar{k}(\rho)=
(||\bar{e}^{n/(n-1)}||^{(n-1)/{n}}_{(N_{r+\ep}),B_\rho}+||\bar{f}||_{(N_{r+\ep}),B_{\rho}})^{1/(n-1)}+||\bar{g}||_{(N_{r+\ep}),B_{\rho}}^{1/n}.
\end{equation}
Let $h$ be any of $\bar{e}^{n/(n-1)},\bar{f},\bar{g}$.
By Lemma ~\ref{shrinkS}, there exists a constant $C$ depending only on $||h||_{(N_r), B_{\rho}}$, $r$, and $n$ such that
$||h||_{(N_{r+\ep}), B_{\rho}}\leq C |\log (\rho )|^{-\ep }$.
Hence, we have
\begin{equation}
\bar{k}
\leq C|\log (\rho )|^{-\gamma }, \label{bark}
\end{equation}
with $\gamma ={\ep}/{n}$.

Applying Theorem ~\ref{Harnackthm1} to $\bu$ in $B_{\rho/3}$, we get
\begin{equation}
M-\mu'=\sup_{B_{\rho /3}}\bu \leq C( \inf_{B_{\rho /3}}\bu +\bar{k})=C(M-M'+\bar{k}),\ \forall \rho \in (0,R] ,\label{M}
\end{equation}
where $M'=M'(\rho )=M(\rho /3)$ and $\mu'=\mu'(\rho )=\mu (\rho /3)$. In this inequality, $C$ depends only on the quantities listed
in this theorem (using the norm $||*||_{(N_{r+\ep}),D}$), and in particular is independent of $R\leq R_0$. In what follows, we fix this $C$,
so, in the same way, we have
\begin{equation}
M'-\mu =\sup_{B_{\rho /3}}\bar{\bu}\leq C(\inf_{B_{\rho /3}}\bar{\bu}+\bar{k})=C(\mu'-\mu +\bar{k}).\label{M'}
\end{equation}
Adding \rq{M} and \rq{M'}, we get
\begin{equation}
M'-\mu'\leq \frac{C-1}{C+1}(M-\mu)+\frac{2C\bar{k}}{C+1}.
\label{M''}
\end{equation}

Let $$\theta =\frac{C-1}{C+1},\ \ \ \tau = \frac{2C}{C+1},\ \ \  \omega (\rho ) = M(\rho) -\mu(\rho).$$
Then \rq{M''} becomes
$$\omega (\rho /3)\leq (\theta\omega (\rho )+\tau \bar{k}(\rho) ).$$
We iterate this relation from $\rho=R$ with successively smaller radii, getting
\begin{equation}
\omega (R3^{-m})\leq \theta^m\omega (R)+\tau \sum_{j=1}^{m}\theta^{j-1} \bar{k}(R3^{j-m}),\ m\in \bN.\label{omega}
\end{equation}
We use this to prove Lemma \ref{omegabound} below.
\noindent
\begin{Lemma} There exists $K>0$ as in Theorem \ref{continuous}
such that if $\rho < R^2$, then
$$\omega( \rho )\leq \frac{K}{|\log(\rho )|^{\gamma}}.
$$
\label{omegabound}
\end{Lemma}
Assuming the lemma for now, by standard arguments we can redefine $u$ at any non-Lebesgue points in $D'$ using limits of averages. Lemma \ref{omegabound} implies these limits converge,
and that the new function is continuous on $D'$. Using compact subsets to cover $D$, we can get continuity on all of $D$. Theorem \ref{continuous} follows; to prove \rq{ContEqn} we apply the lemma with any $\rho$ such that $|x-x_0|<\rho < R^2$.
\hfill $\Box$
\medskip
\noindent
{\it Proof of lemma:}
In what follows, we use $K$ to denote various positive constants that are independent of $u,R$.
To analyze \rq{omega}, we apply \rq{bark},
$$\bar{k}(R3^{j-m}) \leq K|{\log(R3^{j-m})}|^{-\gamma } = K(|\log R|+(m-j)\log 3)^{-\gamma } \leq \frac{K}{(m-j)^\gamma }$$
when $1\le j<m$. For $j=m$, we get $\bar{k}(R) \leq K|{\log R}|^{-\gamma } \le K$.
Then
$$\sum_{1\le j \le m/2} \theta^{j-1} (m-j)^{-\gamma}\leq K m^{-\gamma}$$
since $m-j \ge m/2$ and $\sum_{j=1}^\infty \theta^{j-1}$ converges. Furthermore, using $(m-j)^{-\gamma} \le 1$,
$$\sum_{m/2 < j \le m-1} \theta^{j-1}(m-j)^{-\gamma}\leq K\theta^{m/2}.$$
Applying these to \rq{omega} (and leaving the $j=m$ term to the reader),
\begin{eqnarray*}
\omega(R 3^{-m})& \leq & \theta^m\omega(R)+\tau K  m^{-\gamma}+\tau K\theta^{m/2}.
\end{eqnarray*}
Since
$\omega(R) \le 2||u||_{L^\infty(D'')} <K$
and max$_m \ {m^\gamma} \theta^{m/2} \le K$, we have
$$\omega(R 3^{-m}) \leq \frac{K}{m^\gamma}.$$
Since $R\le R_1 \le 3^{-10}$, we have  $\rho \leq R3^{-10}$.
Fix $m\ge 10$ such that $R3^{-m-1}< \rho \leq R3^{-m}$.
It is not hard to show that $\rho /R>3^{-m-1}$ implies $m\geq (.9/\log 3)|\log (\rho /R)|$. Since $\omega (\rho)$ is a non-decreasing function of $\rho$, we have
$$\omega (\rho )\leq \omega(R 3^{-m})\leq \frac{K}{m^\gamma} \leq \frac{K}{|\log (\rho /R)|^\gamma}.$$
Since $\rho \leq R^2$, we have $|\log (\rho /R)|^{-1}\leq 2|\log (\rho )|^{-1}$, proving the lemma. \hfill $\Box$
\medskip

We conclude this section by comparing our results here to their counterpart in Serrin's paper,\cite{S}. If one supposes that $h\in L^{p+\ep}(B_R)$, with $p\geq 1$, then H\"older's inequality
gives
$$||h||_{p,B_R}\leq CR^{n\ep /(p(p+\ep ))}||h||_{p+\ep, B_R},$$
and thus $||f||_{p,B_R}$ will vanish algebraically quickly as $R$ goes to zero, a stronger conclusion that in our Lemma \ref{shrinkS}. With this stronger estimate, Serrin
was able to prove H\"older continuity for his solution to \rq{SerrinPde}. Inspecting our proof, it is not hard to see that  if we strengthened our hypotheses on coefficient functions $b, d, e, f,$ and $g$ (but not necessarily on $c$)
 so that
$\bar{k}(\rho )=O(\rho ^{\delta})$ for
some $\delta >0$, then we could also prove H\"older continuity. For works that consider the regularity of u in the case
$b=d=e=f=g=0$, see \cite{NU}, \cite{SSSZ}.
\end{section}

\begin{section}{Appendix}
\begin{subsection}{Proof of Lemma \ref{WkTestFn}}

We assume that
\beq
a>0 \mbox{ is constant, with }b^{n/(n-1)}, d, f, c^n \in L\log L^{(n-1)}(D),\mbox{ and } e \in L^{n/(n-1)}(D),
\lab{r=1}
\eeq which is slightly weaker than \rq{a-g}. Assume \rq{Serrin6} and that
\rq{IntABphi} holds for $\phi \in C_0^{\infty}(D)$.
Let $v\in W_0^{1,n}(B)$ with $B = B_R(z)$ and $\overline{B_{R}(z)}\subset D$.
We need to show that
\rq{IntABphi} holds for $v$.

First, we claim that $||u||_{E_n,B}<\infty$.
Let $h >0$ with $B' = B_{(1+h)R}(z)\subset D$.
Let $\eta \in C^{\infty}_0 (B')$ such that $\eta =1$ on $B$ and $|\eta_x|\leq 2/h$. Then, $\eta u \in W^{1,n}_0(B')$,
and by Lemma \ref{MTEn},
$||u||_{E_n,B} \leq ||\eta u ||_{E_n,B'} \leq C||(\eta u)_x||_{n,B'} \leq C(||u||_{n,B'} + ||u_x||_{n,B'})$, proving the claim.

Unless stated otherwise, all integrals and norms below are over $B$.
We will show that the first integral converges, and then show it is zero. By \rq{Serrin6},
\begin{eqnarray*}|\int_{B} v_x \cdot {\cal A}(x,u,u_x) - v {\cal B}(x,u,u_x) \; dx|
&\leq & \int_B |v_x||{\cal A}| + \int_B |v| |{\cal B}|\\
& \leq & a\int_B  |u_x|^{n-1}|v_x| \;dx +\int_B b(x)|u|^{n-1}|v_x| \;dx + \int_B e(x)|v_x| \;dx,\\
& + & \int_B c(x) |u_x|^{n-1}|v| \;dx +\int_B d(x)|u|^{n-1}|v| \;dx + \int_B f(x)|v| \;dx\\
& = & I_a + I_b + I_e + I_c + I_d + I_f.
\end{eqnarray*}
We will prove that
\beq
I_a,\;I_b,\;I_c,\;I_d\;,I_e,\;I_f \leq C||v_x||_n\ .
\lab{Ia-If}
\eeq
Assuming this for the moment, all integrals above must converge.
Also, let $\{ \phi_j \}$ be a sequence of  $C^{\infty}_0 (B)$ functions converging to $v$ in $W_0^{1,n} (B)$.
Applying \rq{Ia-If} with $v$ replaced by  $v-\phi_j$,
\begin{eqnarray*} |\int_{B} v_x \cdot {\cal A}(x,u,u_x) - v{\cal B}(x,u,u_x)\;dx| \le ||(v - \phi_j)_x||_n \rightarrow 0
\end{eqnarray*}
which proves the Lemma.
\vspace{3mm}

Returning to \rq{Ia-If}, the estimates for $I_a,\;I_b,\;I_c,\;I_e,\;I_f $ follow from fairly straightforward arguments similar to those for $I_a,\;I_b,\;I_c,\;I_e $ in the proof of Theorem \rq{bounded}.
We include the details only for $I_d$.
By \cite{Leck}, if $w\in W_0^{1,n}(B')$ with $||w_x||_{n,B'} \leq 1$, then
\beq
\int_{B'} e^{\alpha_n |w|^{n/(n-1)} } \;dx \leq C|B'|,
\lab{noluxMT}
\eeq
where $\alpha_n > 0$ is a fixed constant. Let $\hat{M}(t) = M_1((\alpha_n/2)  t)$ and $\hat{N}(t) = N_1(2t/\alpha_n)$. Then $(\hat{M}, \hat{N})$ is an Orlicz pair of convex functions and, it follows from $d\in L\log^{n-1}L(B')$, that $\int_B \hat{N}(|d|)  <\infty$.

Set $w = \eta u/||(\eta u)_x||_{n,B'}$ so that $||w_x||_{n,B} \leq ||w_x||_{n,B'} = 1$.
Set $\lambda = ||v_x||_n$. Using Young's inequality for $(\hat{M},\hat{N}) $, $\hat{M}(t)\leq Ce^{\alpha_n t^{1/(n-1)}}$, the inequality $|ab|\leq \frac{n-1}{n}a^{n/(n-1)}+\frac{1}{n}b^{n}$, H\"older's inequality, and \rq{noluxMT}:
\begin{eqnarray*}I_d
& = & \left(\int_B d(x)|w^{n-1} \, v| \;dx\right) ||(\eta u)_x||^{n-1}_{n,B'}\\
& \le & C \lambda \int_B |\frac{w^{n-1}
v}{\lambda}|\ |d| \; dx \\
& \leq & C \lambda \left(\int_B \hat{M}\left(\frac{|w|^{n-1}
 v|}{\lambda}\right) dx+ \int_B \hat{N}(|d|) dx \right)\\
 & \leq &  C\lambda\left( \int_B \exp \left(\alpha_n |w|\left(\frac{
|v |}{\lambda}\right )^{1/(n-1)}\right )dx + 1\right) \\
 & \leq & C\lambda \left( \int_B
 e^{\frac{n-1}{n}\alpha_n|w|^{\frac{n}{n-1}}}
 e^{\frac{1}{n}\alpha_n(|v|\lambda^{-1})^{\frac{n}{n-1}}} dx+ 1 \right)\\
 & \leq & C\lambda\left\{ \left(\int_B e^{\alpha_n|w|^{{\frac{n}{n-1}}}}dx\right)^{\frac{n-1}{n}}
 \left(\int_B e^{\alpha_n(\frac{|v|}{\lambda})^{{\frac{n}{n-1}}}}
 dx\right)^{\frac 1 n}  +1\right\}\\
 & \leq & C||v_x||_n. \hspace{4.5in} \hfill  \Box
 \end{eqnarray*}
\end{subsection}

\begin{subsection}{Proof of Lemmas \ref{MTEn}, \ref{MonoConv} and \ref{boundedE}.}\label{MT1}
\noindent
{\it Proof of Lemma \ref{MTEn}:}
To prove parts A and B, we define $t$ by  $e^{-t/n} = |x|/R $.
Then $d|x|/dt = -(R/n)e^{-t/n} = -|x|/n$.
Let $u^{\#}(x) = u^{\#}(|x|)$
be the  symmetric-decreasing rearrangement of $|u|$ on $B_R (0)$; see \cite{LL}.
So, $|| u^{\#}_x||_n \leq || u_x||_n$ and, for any of our Orlicz
norms, $||u|| = ||u^{\#}||$.
Let $\sigma$ be the $(n-1)$ dimensional measure of the unit sphere in $\bR^n$, so that
$|B_R| = \sigma R^n/n$.  Define $\omega(t) = C_nu^{\#}(|x|) $, where $C_n= (\sigma n^{(n-1)})^{1/n}$, and

\begin{equation}\label{mark1}
 \int_0^{\infty} (\omega'(t))^n dt = || u_x^{\#}||_n^n \ .
\end{equation}
Let $0 < r \leq 1$ (the case $r=1$ is used for part B.
We have $$||u||_{E_{nr}} = (C_n)^{-1}||C_n u||_{E_{nr}} = (C_n)^{-1}\inf \{\la : \int_{B_R} e^{ (C_n |u|/\la )^{nr/(n-1)} } -1\;dx  \leq 1  \}.$$ For $\la  > 0$,
\begin{eqnarray}
 I_A: & = & \int_{B_R} e^{(C_n |u|/\la )^{nr/(n-1)}} - 1 \;dx\nonumber \\
 &  =&  \int_{B_R} e^{(C_n u^{\#}/\la )^{nr/(n-1)}} - 1 \;dx,\nonumber \\
& = & (\sigma R^n/n)\int_0^R (e^{(C_n u^{\#}(|x|)/\la )^{nr/(n-1)}} - 1) \;d(|x|^n/R^n),\nonumber \\
& = & |B_R| \int_0^{\infty} (e^{( \omega(t)/\la )^{rn/(n-1)}} - 1)e^{-t} \;dt .\label{mark2}
\end{eqnarray}
Using $\omega(0) = 0$ and H\"older's inequality, we have
$$\omega(t) = \int_0^t \omega'(\tau)\;d\tau \leq (t)^{(n-1)/n}(\int_0^t (\omega'(\tau))^n\;d\tau)^{1/n} \leq t^{(n-1)/n}||u_x||_n.$$
Let $K = (||u_x||_n/\la )^{n/(n-1)}$, so that
\begin{equation}\label{mark3}
I_A \leq |B_R|\int_0^{\infty} (e^{(Kt)^r} - 1)e^{-t}\;dt := |B_R|I_K.
\end{equation}

\noindent
{\it Proof of A:}
Let $r<1$. To bound $||u||_{E_{nr}}$, we will need a good $\la > 0$, based on an upper bound for $I_K =\int_0^{\infty} (e^{(Kt)^r} - 1)e^{-t}\;dt$.
We will use the identities $s! = \Gamma (s+1) = \int_0^{\infty} t^s e^{-t} \; dt$ and $ e^{x} - 1  = \Sigma_{j = 1} \;x^j/j!$. These identities and H\"older's inequality give,
\begin{eqnarray*}
I_K & = & \sum_{j = 1}\; \int_0^{\infty} (Kt)^{jr}\; e^{-t}/j!  \; dt \\
& = & \sum_{j = 1}\;  (K)^{jr} (jr)!/j! \\
& \leq & (\sum_{j = 1}\;  (K)^{jr/(1-r)}/j!)^{1 - r} (\sum_{j = 1}\;  [(jr )!]^{1/r}/j!)^{r} \\
& = &  (e^{K^{r/(1-r)}}-1)^{1 - r} S^r.
\end{eqnarray*}

\vspace{3mm}
We show that $S =  \Sigma_{j = 1}\;  [(jr )!]^{(1/r)}/j!$ converges for $0< r< 1$. By Stirling's formula, $m! \approx (2\pi m)^{1/2} m^m e^{-m} $, so the terms of this series compare to,
$$\frac{[(jr )!]^{1/r}}{j!} \approx \frac{  (2\pi jr)^{1/(2r)} (jr)^j e^{-j}}{(2\pi j)^{1/2} (j)^j\; e^{-j}} = Cj^{(1-r)/(2r)} r^j.$$
(where $C$ may depend on $r$ but not $j$). Now set $C_r$  to be the maximum of $\{ 1, S^{r/(r-1)}\}$. Then

 $$ I_K \le (C_r(e^{K^{r/(1-r)}} -1))^{1 - r} \leq (e^{ C_r K^{r/(1-r)}} -1)^{1-r}.$$

We have $I_A \leq |B_R|(e^{ C_r K^{r/(1-r)}} -1)^{1-r}$ for $K = (||u_x||_n/\la )^{n/(n-1)}$. Set
$$|B_R|(e^{ C_r ( (||u_x||_n/\la )^{n/(n-1)})^{r/(1-r)}} -1)^{1-r} = 1,$$
to determine $\la > 0$. Since $C_r ( (||u_x||_n/\la )^{n/(n-1)})^{r/(1-r)} = \log( 1 + [1/|B_R|]^{1/(1-r)})$, we get
$$ \la  = C^{\delta}_r (\log( 1 + 1/|B_R|^{1/(1-r)}))^{-\delta }|| u_x||_n ,$$ where $\delta = \frac{(n-1)(1-r)}{nr}$. So,
$$||u||_{E_nr} \le \frac{C^{\delta}_r}{C_n} (\log( 1 + [1/|B_R|]^{1/(1-r)}))^{-\delta }|| u_x||_n.$$
\vspace{5mm}

\noindent
{\it Proof of B:}
We begin by arguing exactly as in part A
up to \rq{mark3}, but now let $r=1$.
Define $K$ by $|B_R|K/(1-K) = 1$, which also determines $\la$. Then
$$I_K = \int_0^{\infty} (e^{Kt} - 1)e^{-t}\;dt = 1/(1-K) - 1 = |B_R|^{-1}$$
so that $I_A \le |B_R|I_K = 1$ and $||C_n u||_{E_n} \leq \la   = ||u_x||_n(|B_R| +1)^{(n-1)/n}$, proving part B.
\hfill $\Box$
\vspace{5mm}

\noindent
{\it Proof of Lemma \ref{MonoConv}:}
Let $k_l :=||v||_{E_m,h'}$, and assume for now that $k_0 :=||\bar{u}^q||_{E_m,h'} < \infty$.
Since $F'' \le 0$,
we have $v\leq \bar{u}^q$ so that $k_l\leq k_0$ for all $l$, and $k_l$ increases with $l$.
So, $k_{\infty} :=\lim_{l\to \infty}k_l \leq k_0$.
Let $\ep >0$. By the definition of $||*||_{E_m}$,
$$ I_v := \int_D \exp \left ( (\frac{v(x)}{k_{\infty}+\ep})^{\frac{m}{n-1}}\right ) -1 \ dx \leq 1.$$
By Monotone Convergence, $I_{\bar{u}^q}\leq 1.$ as well,
which implies $k_{\infty}+\ep \geq k_0$ and that $k_{\infty}=k_0$.
It is not hard to see the same proof works if $k_0=\infty$.
\hfill $\Box$
\vspace{5mm}

\noindent
{\it Proof of Lemma \ref{boundedE}:}
Assume $||\bar{u}||_{\infty,R} >J+\epsilon>J$, to get a contradiction.
So, $\bar{u}> J+\epsilon$ on some set $\Omega \subset B_R(0)$ of positive measure.
From the definition of the norm, for each $j$, $\int_\Omega \exp(|\bar{u} / J|^{q_j/(n-1)})-1\ dx \le 1$.
So, $|\Omega| (\exp((1+ \epsilon/ J)^{q_j/(n-1)})-1) \le 1$.
Letting $j\to \infty$, we get a contradiction.
\hfill $\Box$

\end{subsection}
\end{section}

\end{document}